\documentclass{amsart}
\usepackage{tikz,tikz-cd}
\usepackage{url}
\usepackage{graphicx}
\usepackage{color}
\usepackage{times}
\usepackage{cite}
\usepackage{enumerate,latexsym}
\usepackage{amsmath,amssymb}
\usepackage{amsthm}
\usepackage{verbatim}
\usepackage{mathrsfs}
\usepackage{tabularx}

\newtheorem{thm}{Theorem}[section]
\newtheorem*{theorem*}{Theorem}
\newtheorem*{acknowledgement*}{Acknowledgement}

\newtheorem{lem}[thm]{Lemma}
\newtheorem{prop}[thm]{Proposition}
\theoremstyle{definition}
\newtheorem{defn}[thm]{Definition}
\theoremstyle{remark}
\newtheorem{rem}[thm]{Remark}

\numberwithin{equation}{section}

\newcommand{\set}[1]{\left\{#1\right\}}
\newcommand{\Real}{\mathbb R}

\newcommand{\func}[1]{\ensuremath{\mathop{\mathrm{#1}}} }

\newcommand{\Div}[0]{\func{div}}

\newcommand{\spt}[0]{\func{spt}}

\newcommand{\cC}[0]{\mathcal{C}}

\title[Density of Minimal Cones]{Lower Bounds on Density for Topologically Nontrivial Minimal Cones up to Dimension Six}

\author{Jacob Bernstein}
\address{Department of Mathematics, Johns Hopkins University, 3400 N. Charles Street, Baltimore, MD 21218}
\email{jberns15@jhu.edu}

\author{Lu Wang}
\address{Department of Mathematics, Yale University, 219 Prospect Street, New Haven, CT 06511}
\email{lu.wang@yale.edu}

\thanks{The first-named author was partially supported by the NSF Grant DMS-1904674 and DMS-2203132. The second-named author was partially supported by the NSF Grant  DMS-2146997.}

\begin{document}

\begin{abstract}
We prove lower bounds on the density of regular minimal cones of dimension less than seven provided the complements of the cones are topologically nontrivial.
\end{abstract}

\maketitle

\section{Introduction} \label{IntroSec}
In this paper, we prove nearly sharp lower bounds on the density of certain minimal cones of dimension less than $7$.   Recall, a \emph{regular minimal cone} in $\mathbb{R}^{n+1}$ is a cone,  $\mathcal{C}$,  with vertex at the origin $\mathbf{0}$ such that $\mathcal{C}\setminus\{\mathbf{0}\}$ is a non-empty smooth minimal hypersurface -- i.e., its mean curvature is equal to $0$ at all points.  When $n\geq 2$,  a regular cone is minimal if and only if its associated varifold is stationary for area but this is not necessarily the case when $n=1$. The \emph{density},  $\Theta(\mathcal{C})$,  of $\mathcal{C}$ at $\mathbf{0}$ is defined to be
\begin{equation} \label{DenConeEqn}
    \Theta(\mathcal{C})=\frac{\mathcal{H}^{n}(\mathcal{C}\cap B^{n+1}_1(\mathbf{0}))}{\omega_n} =\frac{\mathcal{H}^n(\mathcal{C}\cap B^{n+1}_{R}(\mathbf{0}))}{\omega_n R^n}, R>0,
\end{equation}
where $B^{n+1}_R(\mathbf{0})$ is the open ball in $\mathbb{R}^{n+1}$ centered at $\mathbf{0}$ with radius $R$, $\omega_n=|B_1^n(\mathbf{0})|$ is the volume of the unit $n$-ball in $\Real^n$, and $\mathcal{H}^n$ is the $n$-dimensional Hausdorff measure. The upper semi-continuity of density for stationary varifolds implies that $\Theta(\mathcal{C})\geq 1$ and standard dimension reduction arguments ensure equality occurs only when $\mathcal{C}$ is a hyperplane.   Allard's regularity theorem \cite{AllReg} implies that there exist constants $\epsilon(n)>0$ so that if $\cC$ is a non-flat regular minimal cone, then  $\Theta(\cC)\geq 1+\epsilon(n)$.   In \cite{ChengLiYau}, Cheng--Li--Yau  give an explicit, but very rough, lower bound for $\epsilon(n)$ -- see Appendix \ref{Explicit}.
 
Following Colding--Minicozzi \cite{CMGenMCF}, the \emph{entropy} of a hypersurface $\Sigma\subset\mathbb{R}^{n+1}$ is
\begin{equation} \label{CMEntEqn}
    \lambda(\Sigma)=\sup_{\mathbf{x}_0\in\mathbb{R}^{n+1}, t_0>0} (4\pi t_0)^{-\frac{n}{2}} \int_\Sigma e^{-\frac{|\mathbf{x}-\mathbf{x}_0|^2}{4t_0}} \, d\mathcal{H}^n.
\end{equation}
Because a regular stationary cone, $\mathcal{C}$, may be thought of as an eternal weak mean curvature flow, Huisken's monotonicity formula \cite{HuiAsympMCF, IlmSingMCF} implies that $\lambda(\mathcal{C})=\Theta(\mathcal{C})$.  Likewise, the entropy of a round $k$-sphere, $\mathbb{S}^k\subset \mathbb{R}^{k+1}$ equals the Gaussian density of the self-similarly shrinking $\mathbb{S}^k$ at its center \cite{CMGenMCF}.  Thus, by Stone's computation \cite[Appendix A]{StoEntSphere},
\begin{equation} \label{EntSphereEqn}
    2=\lambda(\mathbb{S}^0)>\lambda(\mathbb{S}^1)>\lambda(\mathbb{S}^2)>\lambda(\mathbb{S}^3)>\cdots\to\sqrt{2}.
\end{equation}

In  \cite[Theorem $\text{1}^*$]{IWDenMinCone}, Ilmanen--White used mean curvature flow and the existence of the Hardt--Simon foliation \cite{HSFoliation} to show that if $\mathcal{C}$ is a regular area-minimizing cone that is topologically non-trivial, then $\Theta(\mathcal{C})\geq \lambda(\mathbb{S}^{n-1})>\sqrt{2}$.  There are no non-flat area-minimizing cones when $n\leq 6$ and so their theorem does not apply in these dimensions. However, using a different argument inspired by \cite{BWTopolUnique} we obtain the same lower bound for any topologically non-trivial regular minimal cones in precisely this range of dimensions.
\begin{thm} \label{NonContractThm}
    For $ n\leq 6$, let $\mathcal{C}$ be a regular minimal cone in $\mathbb{R}^{n+1}$. If at least one of the components of $\mathbb{R}^{n+1}\setminus\mathcal{C}$ is not contractible, then 
    \[
     \Theta(\mathcal{C})\geq \lambda(\mathbb{S}^{n-1})>\sqrt{2}.
    \]
\end{thm}
This partially answers a question raised in \cite[\S 5, Problem 1]{IWDenMinCone}.  Using the same method together with work of White \cite{WhiTopolMCMCF} we also show that under  stronger topological hypotheses one obtains a better lower bound.
\begin{thm} \label{NonHomolBallThm}
    For $ 2\leq n\leq 6$, let $\mathcal{C}$ be a regular minimal cone in $\mathbb{R}^{n+1}$. If at least one of the components of $\mathbb{R}^{n+1}\setminus \mathcal{C}$ is not a homology ball, then 
    \[
     \Theta(\mathcal{C}) \geq \lambda(\mathbb{S}^{n-2})>\sqrt{2}.
    \]
\end{thm}
In very low dimensions these bounds are trivial. Indeed, a regular minimal cone in $\mathbb{R}^2$ is the union of $\ell$ of rays based at $\mathbf{0}$ -- though the associated varifold is not stationary unless a balancing condition is satisfied. In particular, every component of the complement of a regular minimal cone is contractible and so Theorem \ref{NonContractThm} is vacuous. We observe that the entropy of such a stationary cone is $\frac{\ell}{2}$, but can be higher for general regular minimal cones. It is not hard to see that the lowest density of a non-trivial stationary cone in $\mathbb{R}^2$ is $\frac{3}{2}$. Moreover, $\ell$ is even if and only if the associated varifold is cyclic mod 2; the condition of being cyclic mod 2 is automatically satisfied by any regular minimal cone in higher dimensions. Hence, the lowest entropy of a non-trivial stationary and cyclic mod 2 cone in $\Real^2$ is $\lambda(\mathbb{S}^0)=2$. Likewise, in $\mathbb{R}^3$, the only regular minimal cones are planes because great circles are the only closed geodesics in $\mathbb{S}^2$. Thus, Theorems \ref{NonContractThm} and \ref{NonHomolBallThm} are again both vacuous.  Within the larger class of non-trivial stationary cones in $\Real^3$ one readily sees that the density is bounded below by $\frac{3}{2}$ in general and by $2$ when the cones are required to be cyclic mod 2.  We remark that the cone over the edges of a regular tetrahedron is a stationary cone with density lying in $(\frac{3}{2}, 2)$.

In the first non-trivial dimension, $\Real^4$, the classification of surfaces and Alexander's theorem ensure that the hypotheses of Theorems \ref{NonContractThm} and \ref{NonHomolBallThm} are both equivalent to the hypothesis that the link of the cone has positive genus and so, in this dimension, Theorem \ref{NonHomolBallThm} implies Theorem \ref{NonContractThm}.    When $n\geq 4$ there exist homology balls that are not contractible and so the hypotheses of  Theorem \ref{NonContractThm} are genuinely weaker than those of Theorem \ref{NonHomolBallThm}. In \cite{Hsiang1, Hsiang2,HsiangSterling, HsiangTomter} there are many examples of non-trivial regular minimal cones whose links are topological spheres.  The authors' are unaware an example of a regular minimal cone whose link is a homology ball but is not a homotopy ball.

Ilmanen--White \cite[Theorem 2]{IWDenMinCone} also proved that given a regular area-minimizing cone, if one of the components of the complement of the cone has nontrivial $k$-th homotopy group, then the density of the cone at $\mathbf{0}$ is greater than or equal to $\lambda(\mathbb{S}^k)$. By Theorem \ref{NonHomolBallThm}, this result is also true for regular minimal cones in $\mathbb{R}^4$, but is stronger than what we are able to show when $4\leq n \leq 6$ and it seems that extending this result to non-area minimizing cones remains an open problem. We remark that, in his thesis \cite{ZhuThesis}, Zhu was able to obtain non-trivial bounds under the weakest possible topological hypotheses.  Namely, it is a simple consequence of \cite[Corollary 2.2]{ZhuThesis} that if $\cC$ is a $n$-dimensional regular minimal cone whose link is not isotopic to the standard sphere, then $\Theta(\cC)\geq \frac{\lambda(\mathbb{S}^{n-2})}{\lambda(\mathbb{S}^{n-1})}>1$.  When $n\geq 4$, the topological hypotheses of Theorem \ref{NonHomolBallThm} are only enough to ensure the link is not a homology sphere. However, the weaker topological hypothesis in Zhu's result results in a worse lower bound on density -- see Appendix \ref{Explicit}.  

The density bounds of Theorems \ref{NonContractThm} and \ref{NonHomolBallThm} are probably not sharp for cones of a given dimension but cannot be improved by much. For positive integers $l$ and $m$, let 
$$
\mathcal{C}_{l,m}=\set{(\mathbf{y}, \mathbf{z})\in \Real^{m+1}\times\Real^{l+1}: m|\mathbf{y}|^2=l|\mathbf{z}|^2}\subset \Real^{m+l+2}
$$
be the family of generalized Simons cones. One readily computes -- see Appendix \ref{Explicit} -- that $\Theta(\cC_{1,1})=\pi/2\approx 1.57$.  In fact,  by Marques--Neves' proof of the Willmore conjecture, i.e.,  \cite[Theorem B]{MNWillmore}, and Almgren's theorem \cite{AlmMinSurfSphere} this is the sharp lower bound for any non-flat regular minimal cones in $\Real^4$.  By comparison, the bound coming from Theorem \ref{NonHomolBallThm} is $\lambda(\mathbb{S}^1)=\sqrt{2\pi}/e\approx 1.52$.  Thus,  the lower bound from Theorem \ref{NonHomolBallThm} is only about $3 \%$ lower than the sharp bound. Likewise, for $4\leq n \leq 6$ and $q=\frac{n-1}{2}$, the cone $\mathcal{C}_{\lfloor q \rfloor, \lceil q \rceil}$ maximizes density among the $\mathcal{C}_{l,m}$ with $l+m=n-1$, and its complement contains a component that is not a homology ball and so also not contractible. It is readily checked that the density of $\mathcal{C}_{\lfloor q \rfloor, \lceil q \rceil}$ at $\mathbf{0}$ is about $3\%$--$5\%$ higher than the lower bounds from Theorems \ref{NonContractThm} and \ref{NonHomolBallThm}.  Indeed, since $\Theta(\mathcal{C}_{l,m})\to \lambda(\mathbb{S}^l)$ as $m\to\infty$, the bounds $\lambda(\mathbb{S}^{n-1})$ and $\lambda(\mathbb{S}^{n-2})$ of Theorems \ref{NonContractThm} and \ref{NonHomolBallThm} may be the best possible that are independent of dimension. Similarly, $\Theta(\mathcal{C}_{l,m})\to \sqrt{2}$ as $l,m\to \infty$. According to \cite{IWDenMinCone},  $\sqrt{2}$ was conjectured by B. Solomon to be the optimal lower bound on the density of non-trivial regular area-minimizing cones.  In \cite[pg. 288]{YauBook}, S. T. Yau asked the more ambitious question of whether appropriate $\mathcal{C}_{l,m}$ minimize area among non-totally geodesic minimal hypersurfaces in the sphere --  in \cite{YauBook} this question is attributed to  B. Solomon as a conjecture -- see also \cite[Section 5]{IWDenMinCone}.  In \cite{RotAreaLB}, this stronger question is answered in the affirmative among highly symmetric minimal hypersurfaces.

To prove Theorems \ref{NonContractThm} and \ref{NonHomolBallThm}, we use properties of self-expanders -- i.e., solutions to \eqref{ExpanderEqn}.  This is because the existence of the Hardt--Simon foliation \cite{HSFoliation}, one of the key ingredients in the approach of \cite{IWDenMinCone}, requires the cone to be area-minimizing and not just minimal.  Instead, for a regular minimal cone that is not area-minimizing, it follows from \cite{BWTopolUnique} and \cite{DinMCExpander}, that there are two self-expanders asymptotic to the cone with the property that any other self-expander asymptotic to the cone is trapped between them. Furthermore, the complement of each of these self-expanders is the union of two components, one star-shaped relative to $\mathbf{0}$ and the other homotopy equivalent to a component of complement of the link.  By combining ideas from \cite{BWTopolUnique} and \cite{BCWMorseFlow}, we can then show, in low dimensions and in a  certain generic sense, that there is a finite collection of gradient flows for the expander functional (see \eqref{ExpanderEnergyEqn}) whose union essentially connects the two self-expanders above -- here the dimension restriction is related to regularity properties of minimizing hypersurfaces that simplify things but do not seem to be essential to the argument. These flows evolve in a monotone manner and exhibit good regularity properties. This allows us to use work of White \cite{WhiTopolMCMCF} to establish a relationship between the topologies of the self-expanders and the entropy of the cone. To complete the proof, we show it is possible to reduce to the generic situation by a suitable perturbation of the two self-expanders and their asymptotic cones. 

We point out that while Ilmanen--White's argument uses properties special to area-minimizing cones, our argument uses properties that are special to regular minimal cones that are not area-minimizing.  Another key difference between our argument and that of \cite{IWDenMinCone}, and one that explains why we cannot prove as strong a result as \cite[Theorem 2]{IWDenMinCone}, is that in their paper the authors are able to apply \cite{WhiTopolMCMCF} to a single monotone flow while in our paper we apply it to a sequence of monotone flows whose directions alternate.  

This paper is organized as follows. In Section \ref{MonotoneExpanderFlowSec}, we give necessary background on monotone gradient flow for the expander functional. In Section \ref{TopolChangeSec}, we use the results of White \cite{WhiTopolMCMCF} to study the relationship between topologies of self-expanders asymptotic to a given cone of low entropy. In Section \ref{DenMinimalConeSec}, we prove topological properties for self-expanders asymptotic to a regular minimal cone of low entropy. The main results about densities of minimal cones then follow from this.

\subsection*{Notation and conventions} 
Throughout the paper, $B^{k}_r(p)$ and $\bar{B}^{k}_r(p)$ are respectively the open and closed ball of center $p$ and radius $r$ in $\mathbb{R}^{k}$. We omit the superscript $k$ when its value is clear from context. We also omit the center when it is the origin. Denote by $\mathrm{int}(A), \mathrm{cl}(A)$, and $\partial A$ respectively the interior, closure, and boundary of a set $A\subseteq\mathbb{R}^{k}$. 

Unless otherwise specified, the vertex of a cone will always be assumed to be the origin. A cone $\mathcal{C}\subset\mathbb{R}^{n+1}$ is \emph{$C^\gamma$-regular} if the link $\mathcal{L}(\mathcal{C})$ of $\mathcal{C}$ is an $(n-1)$-dimensional embedded $C^\gamma$-submanifold of $\mathbb{S}^n$. Here, when $\gamma$ is not an integer, $C^\gamma$ is understood as the usual H\"{o}lder regularity $C^{\lfloor\gamma\rfloor,\{\gamma\}}$.

A hypersurface $\Sigma\subset\mathbb{R}^{n+1}$ is  \emph{$C^\gamma$-asymptotically conical} if there is a $C^\gamma$-regular cone $\mathcal{C}\subset\mathbb{R}^{n+1}$ such that $\lim_{\rho\to 0^+} \rho\Sigma=\mathcal{C}$ in $C^\gamma_{loc}(\mathbb{R}^{n+1}\setminus\{\mathbf{0}\})$, i.e., there is a smooth hypersurface $\Gamma\subset\mathbb{R}^{n+1}\setminus\{\mathbf{0}\}$ so that in each annulus $B_R\setminus \bar{B}_{R^{-1}}$, for sufficiently small $\rho>0$ the $\rho\Sigma$ and $\mathcal{C}$ can be written as the normal graphs of functions $u_i$ and $u$, respectively, over $\Gamma$ so $u_i\to u$ in the $C^\gamma$ topology. In this case $\mathcal{C}$ is called the \emph{asymptotic cone} of $\Sigma$ and is denoted by $\mathcal{C}(\Sigma)$.

\section{Monotone expander flows} \label{MonotoneExpanderFlowSec}
In this section, we give background on monotone expander flows that are asymptotic to a cone.  These play an important technical role in the proof of the main results.

A hypersurface $\Sigma\subset\mathbb{R}^{n+1}$ is a \emph{self-expander} if it satisfies the equation
\begin{equation} \label{ExpanderEqn}
    \mathbf{H}_{\Sigma}-\frac{\mathbf{x}^\perp}{2}=\mathbf{0}
\end{equation}
where $\mathbf{x}$ is the position vector, the superscript $\perp$ denotes the projection to the unit normal $\mathbf{n}_{\Sigma}$ of $\Sigma$, and $\mathbf{H}_{\Sigma}$ is the mean curvature given by
\[
\mathbf{H}_{\Sigma}=-H_{\Sigma}\mathbf{n}_{\Sigma}=-{\Div}_{\Sigma} (\mathbf{n}_{\Sigma})\mathbf{n}_{\Sigma}.
\]
A hypersurface $\Sigma$ is a self-expander if and only if the family of homothetic hypersurfaces, $\{\Sigma_t\}_{t>0}=\{\sqrt{t}\Sigma\}_{t>0}$, is a \emph{mean curvature flow}, that is, a solution to the equation
\begin{equation} \label{MCFEqn}
    \left(\frac{\partial\mathbf{x}}{\partial t}\right)^\perp=\mathbf{H}_{\Sigma_t}.
\end{equation}
Self-expanders model the behavior of a mean curvature flow when it emerges from a conical singularity (see \cite{ACINonunique} and \cite{BWRelEnt}), so it is natural to study self-expanders $\Sigma$ asymptotic to cones $\mathcal{C}$ in the sense that $\lim_{t\to 0}\mathcal{H}^n\llcorner\sqrt{t}\Sigma=\mathcal{H}^n\llcorner\mathcal{C}$. By \cite[Proposition 3.3]{BWProper}, if the cone $\mathcal{C}$ is $C^\gamma$-regular for some $\gamma\geq 2$, then $\Sigma$ has quadratic curvature decay and is $C^{\gamma^\prime}$-asymptotic to $\mathcal{C}$ for any $\gamma^\prime\in (0,\gamma)$.

Variationally, self-expanders are critical points for the functional
\begin{equation} \label{ExpanderEnergyEqn}
    E(\Sigma)=\int_{\Sigma} e^{\frac{|\mathbf{x}|^2}{4}} \, d\mathcal{H}^n.
\end{equation}
The associated negative gradient flow is then called an \emph{expander flow}, that is, a family $\{\Sigma_t\}$ of hypersurfaces in $\mathbb{R}^{n+1}$ satisfying the equation 
\begin{equation} \label{ExpanderFlowEqn}
    \left(\frac{\partial \mathbf{x}}{\partial t}\right)^\perp=\mathbf{H}_{\Sigma_t}-\frac{\mathbf{x}^\perp}{2}.
\end{equation} 
In general, an expander flow may become singular in finite time.  However,  various notions of weak solutions to \eqref{ExpanderFlowEqn} are at our disposal which allow us to continue the flow through singularities. For the purposes of this paper we mention two of them: \emph{expander weak flow} of closed sets in $\mathbb{R}^{n+1}$ -- see \cite[\S 11]{HWAvoidance} for vector field $-\frac{\mathbf{x}}{2}$; \emph{expander Brakke flow} of Radon measures associated to $n$-dimensional varifolds in $\mathbb{R}^{n+1}$  -- see \cite[\S 13]{HWAvoidance}. We omit the precise definitions for these weak flows because they are not needed in what follows.

Of particular interest is the following special class of weak expander flows whose existence was shown in \cite{BCWMorseFlow}:
\begin{defn} \label{MonotoneExpanderFlowDefn}
Given $T\in \mathbb{R}$ and a $C^3$-regular cone $\mathcal{C}\subset\mathbb{R}^{n+1}$, a \emph{strongly regular strictly monotone expander weak flow asymptotic to $\mathcal{C}$ with starting time $T$} is a family $\mathcal{S}=\{\Omega_t\}_{t\geq T}$ of closed sets in $\mathbb{R}^{n+1}$ with $M_t=\partial\Omega_t$ satisfying:
\begin{enumerate}
    \item The spacetime track of $\mathcal{S}$, $\bigcup_{t\geq T}\Omega_t\times\{t\}$ is a expander weak flow with starting time $T$;
    \item $\Omega_{t_2}\subseteq\mathrm{int}(\Omega_{t_1})$ for $t_2>t_1\geq T$;
    \item Given $\epsilon>0$ there is a radius $R_0>1$ so that for $t\in [T,\infty)$ there is a $C^2$ function $u(\cdot, t)\colon\mathcal{C}\setminus B_{R_0}\to \mathbb{R}$ satisfying
    \[
    \sup_{p\in\mathcal{C}\setminus B_{R_0}} \sum_{i=0}^2 |\mathbf{x}(p)|^{i-1}|\nabla_\mathcal{C}^iu(p,t)|\leq \epsilon
    \]
    and
    \[
    M_t\setminus B_{2R_0}\subset\set{\mathbf{x}(p)+u(p,t)\mathbf{n}_\mathcal{C}(p)\colon p\in\mathcal{C}\setminus B_{R_0}}\subset M_t;
    \]
    \item For $[a,b]\subset [T,\infty)$, $\{M_t\}_{t\in [a,b]}$ is a partition of $\Omega_a\setminus\mathrm{int}(\Omega_b)$; 
    \item $\{\mathcal{H}^n\llcorner M_t\}_{t\geq T}$ is a unit-regular\footnote{We say a flow is \emph{unit-regular} if near every spacetime point of Gaussian density $1$, the flow is regular in a ball in spacetime -- see \cite{WhiLocReg}.} expander Brakke flow;
    \item If $\mathcal{S}_i$ is a blow-up sequence to $\mathcal{S}$ at a point $X_0\in\mathbb{R}^{n+1}\times (T,\infty)$\footnote{The $\mathcal{S}_i$ are obtained by translating $\mathcal{S}$ by $-X_i$ and then parabolically dilating by $\rho_i$, i.e., $(\mathbf{x},t)\mapsto (\rho_i\mathbf{x},\rho_i^2 t$) for $X_i\to X_0$ in spacetime and $\rho_i\to\infty$.}  that converges to a limit flow $\mathcal{S}'=\{\Omega'_t\}_{t\in\mathbb{R}}$, then $\Omega'_t$ is convex for each $t$, and there is a $T^\prime\in [-\infty,\infty]$ so that 
    \begin{enumerate}
        \item $\mathrm{int}(\Omega_t^\prime)\neq\emptyset$ if $t<T^\prime$, while $\mathrm{int}(\Omega_t')=\emptyset$ if $t=T'$;
        \item The spacetime tracks of the $\mathcal{S}_i$ converge, in $C^\infty_{loc}(\mathbb{R}^{n+1}\times (-\infty,T^\prime))$, to the spacetime track of $\mathcal{S}'$, and $\{\partial\Omega'_t\}_{t<T'}$ is a smooth mean curvature flow;
        \item $\Omega_t^\prime=\emptyset$ for $t>T^\prime$.
    \end{enumerate}
    Furthermore, if $\mathcal{S}'$ is a tangent flow, i.e., the blow-ups of $\mathcal{S}$ are all centered at $X_0$, then $\{\partial\Omega'_t\}_{t\in\mathbb{R}}$ is either a static $\mathbb{R}^n$ or a self-similarly shrinking $\mathbb{S}^\ell\times\mathbb{R}^{n-\ell}$ for some $1\leq \ell\leq n$.
\end{enumerate}
\end{defn}

We will need the following result of \cite{BCWMorseFlow}.

\begin{prop} \label{MonotoneExpanderFlowProp}
    For $n\leq 6$, let $\Omega$ be a closed subset of $\mathbb{R}^{n+1}$ such that $\partial\Omega$ is a $C^3$-asymptotically conical hypersurface with asymptotic cone $\mathcal{C}$. Assume that $\Omega$ is strictly expander mean convex, i.e.,
    \[
    2H_{\partial\Omega}(p)+\mathbf{x}(p)\cdot\mathbf{n}_{\partial\Omega}(p)>0 \mbox{ for $p\in\partial\Omega$} 
    \]
    where $\mathbf{n}_{\partial\Omega}$ is the outward unit normal to $\Omega$. Then there exists a strongly regular strictly monotone expander weak flow $\mathcal{S}=\{\Omega_t\}_{t\geq 0}$ asymptotic to $\mathcal{C}$ which starts with $\Omega_0=\Omega$, and such that the $\Omega_t$ converge, in $C^\infty_{loc}(\mathbb{R}^{n+1})$, to $\Omega'=\bigcap_{t\geq 0}\Omega_t$ and $\partial\Omega'$ is a stable self-expander asymptotic to $\mathcal{C}$. 
        
    In addition, if there is a closed set $\Omega''\subseteq\Omega$ such that $\partial\Omega''$ is a smooth self-expander asymptotic to $\mathcal{C}$, then one can ensure that $\Omega''\subseteq\Omega'$.
\end{prop}

\section{Topological change of monotone expander flows and entropy} \label{TopolChangeSec}
In this section, we observe how the results of White \cite{WhiTopolMCMCF} can be used to understand the relationship between the topologies of asymptotically conical self-expanders when the asymptotic cone has sufficiently small entropy.

Throughout this section $\mathcal{C}\subset\mathbb{R}^{n+1}$ will be a cone of at least $C^3$ regularity. Denote its link by $\mathcal{L}(\mathcal{C})$. It is always possible to find an open subset $\omega_+\subset \mathbb{S}^n\setminus \mathcal{L}(\cC)$ with boundary $\mathcal{L}(\cC)$ -- see  \cite[Section 4]{BWTopolUnique} where the pair $(\omega_+, \mathcal{L}(\mathcal{C}))$ is called a boundary link. There is then a unique choice of unit normal,  $\nu$, on $\mathcal{L}(\mathcal{C})\subset\mathbb{S}^n$ so that $\nu$ is the outward normal to $\omega_+$.  We remark that $-\nu$ is the outward normal to $\omega_-=\mathbb{S}^n\setminus \mathrm{cl}(\omega_+)$ and $(\omega_-, \mathcal{L}(\cC))$ is the only other boundary link associated to $\mathcal{L}(\mathcal{C})$. 

Given a $C^3$-asymptotically conical hypersurface $\Sigma\subset\mathbb{R}^{n+1}$ with $\mathcal{C}(\Sigma)=\mathcal{C}$, define $\Omega_+(\Sigma)\subset\mathbb{R}^{n+1}$ to be the closed set with boundary $\Sigma$ such that as $\rho\to 0^+$ the $\rho\Omega_+(\Sigma)\cap \mathbb{S}^n$ converge as closed sets to $\mathrm{cl}(\omega_+)$. In this case one may orient $\Sigma$ so that the outward normal points out of $\Omega_+(\Sigma)$ and this choice is compatible with the choice of $\nu$ in the obvious manner.  Likewise, let $\Omega_-(\Sigma)=\mathbb{R}^{n+1}\setminus\mathrm{int}(\Omega_+(\Sigma))$. 

For two $C^3$-asymptotically conical hypersurfaces $\Sigma_0,\Sigma_1\subset\mathbb{R}^{n+1}$ with the same asymptotic cone $\mathcal{C}$, we say $\Sigma_0\preceq\Sigma_1$ if $\Omega_+(\Sigma_1)\subseteq\Omega_+(\Sigma_0)$. Let $\mathcal{E}(\mathcal{C})$ be the set of self-expanders asymptotic to $\mathcal{C}$, and thus $(\mathcal{E}(\mathcal{C}),\preceq)$ is a partially ordered set. Observe that by swapping $\omega_+$ with $\omega_-$, which is the same as swapping $\nu$ with $-\nu$, one reverses the role of $\Omega_+(\Sigma)$ and $\Omega_-(\Sigma)$ and reverses the partial order.

We first study the case of a \emph{nondegenerate cone}, i.e., a cone, $\cC$ for which there are no nontrivial Jacobi fields that fix infinity on any element of $\mathcal{E}(\mathcal{C})$.  For such $\mathcal{C}$, all stable self-expanders in $\mathcal{E}(\cC)$ are strictly stable. 

\begin{thm} \label{NondegExpanderTopolThm}
    For $2\leq n\leq 6$, let $\mathcal{C}$ be a nondegenerate $C^4$-regular cone in $\mathbb{R}^{n+1}$ such that $\mathcal{L}(\mathcal{C})$ is connected and, for some $m\in [1, n-1]$,
    \[
    \lambda(\mathcal{C})<\lambda(\mathbb{S}^m\times\mathbb{R}^{n-m}).
    \]
    If $\Gamma_+$ and $\Gamma_-$ are two elements of $\mathcal{E}(\mathcal{C})$ with $\Gamma_-\preceq\Gamma_+$, then the inclusions $\mathfrak{i}^\pm \colon \Omega_\pm(\Gamma_\pm) \to \Omega_\pm(\Gamma_\mp)$ induce homomorphisms 
    \[
    \mathfrak{i}^\pm_* \colon H_k(\Omega_\pm(\Gamma_\pm)) \to H_k(\Omega_\pm(\Gamma_\mp))
    \]
    such that 
    \begin{enumerate}
        \item When $n-m\leq k\leq m$, the maps $\mathfrak{i}^\pm_*$ are bijective ;
        \item When $n-m-1=k\leq m$, the maps $\mathfrak{i}^\pm_*$ are injective;
        \item When $n-m\leq k=m+1$, the maps $\mathfrak{i}^\pm_*$ are surjective.
    \end{enumerate}
\end{thm}

\begin{rem} \label{MaxCylinderRem}
It is an immediate consequence of the main result of \cite{BWTopolUnique} that, for $2\leq n\leq 6$, if $\mathcal{C}$ is any $C^3$-regular cone in $\mathbb{R}^{n+1}$ with $\lambda(\mathcal{C})<\lambda(\mathbb{S}^{n-1}\times\mathbb{R})$, then given $\Gamma_+,\Gamma_-\in\mathcal{E}(\mathcal{C})$ with $\Gamma_-\preceq\Gamma_+$, the inclusions $\mathfrak{i}^\pm\colon \Omega_\pm(\Gamma_\pm)\to\Omega_\pm(\Gamma_\mp)$ are homotopy equivalences and so $\mathfrak{i}^\pm_*$ are all isomorphisms.
\end{rem}

To prove Theorem \ref{NondegExpanderTopolThm} we need several auxiliary lemmas/propositions. If  $\mathcal{M}=\{\mu_t\}_{t\geq T}$ is a family of Radon measures on $\mathbb{R}^{n+1}$ and $X_0=(\mathbf{x}_0,t_0)$ is a spacetime point, then the \emph{Gaussian density of $\mathcal{M}$ at $X_0$},  denoted $\theta(\mathcal{M},X_0)$, is defined to be 
\begin{equation} \label{GaussDenEqn}
   \theta(\mathcal{M},X_0)=\lim_{t\to t_0^-} \int (4\pi (t_0-t))^{-\frac{n}{2}} e^{-\frac{|\mathbf{x}-\mathbf{x}_0|^2}{4(t_0-t)}} \, d\mu_t 
\end{equation}
whenever the limit exists and is finite. Otherwise, we set $\theta(\mathcal{M},X_0)=\infty$.

\begin{lem} \label{GaussDenLem}
	Given $T\in \mathbb{R}$ and a $C^3$-regular cone $\mathcal{C}\subset\mathbb{R}^{n+1}$, let $\mathcal{S}=\{\Omega_t\}_{t\geq T}$ be a strongly regular strictly monotone expander weak flow asymptotic to $\mathcal{C}$ with starting time $T$. If $M_t=\partial\Omega_t$ and for some $m\in [1,n-1]$
	\[
	\lambda(M_T)<\lambda(\mathbb{S}^m\times\mathbb{R}^{n-m})
	\]
then for $\mathcal{M}=\{\mathcal{H}^n\llcorner M_t\}_{t\geq T}$ and every point $X_0\in \mathbb{R}^{n+1}\times (T,\infty)$
	\[
	\theta(\mathcal{M},X_0)\leq \lambda(\mathbb{S}^{m+1}\times\mathbb{R}^{n-m-1}).
	\]
\end{lem}

\begin{proof}
    First, without loss of generality we may assume $T=0$. Observe that if $N_\tau=e^{t/2} M_t$ with $\tau=e^t$, then $\lambda(N_\tau)=\lambda(M_{\log\tau})$ for $\tau\geq 1$ and $\{\mathcal{H}^n\llcorner N_\tau\}_{\tau\geq 1}$ is an integral Brakke flow. Thus, by the Huisken monotonicity (see \cite{HuiAsympMCF} and \cite{IlmSingMCF}),  one has that $\lambda(M_t)$ is decreasing in $t$. In particular, $\lambda(M_t)\leq\lambda(M_0)$ for all $t\geq 0$. By a variant of Huisken's monotonicity formula (see \cite[\S 11]{WhiStrata}), as our hypotheses ensure that every tangent flow to $\mathcal{M}$ at a point $X_0\in\mathbb{R}^{n+1}\times (0,\infty)$ is null, a static $\mathbb{R}^n$ or a self-similarly shrinking $\mathbb{S}^\ell\times\mathbb{R}^{n-\ell}$, one has either $\theta(\mathcal{M},X_0)=0$, $\theta(\mathcal{M},X_0)=1$, or
    \[
    \lambda(\mathbb{S}^\ell\times\mathbb{R}^{n-\ell})=\theta(\mathcal{M},X_0)\leq\lambda(M_0)<\lambda(\mathbb{S}^m\times\mathbb{R}^{n-m}).
    \]
    Combined with \eqref{EntSphereEqn}, it follows that $\theta(\mathcal{M},X_0)\leq\lambda(\mathbb{S}^{m+1}\times\mathbb{R}^{n-m-1})$.
\end{proof}

Let $\Omega$ be a closed subset of $\mathbb{R}^{n+1}$ with smooth boundary, and let $K$ be a closed subset of $\Omega$. Following White \cite[Definition 5.1]{WhiTopolMCMCF}, a point $p\in K$ is a \emph{regular point} of $K$ provided that either $p$ is an interior point of $K$;  or $p\in \Omega\setminus\partial\Omega$ and $\Omega$ has a neighborhood, $U$, of $p$ such that $K\cap U$ is diffeomorphic to a closed half-space in $\mathbb{R}^{n+1}$;  or $p\in\partial\Omega$ and $\Omega$ has a neighborhood $U$ of $p$ for which there is a diffeomorphism that maps $U$ onto the closed half-space $\{x_1\geq 0\}\subset\mathbb{R}^{n+1}$ and that maps $K\cap U$ onto $\{x_1\geq 0, x_2\geq 0\}$. Points in $K$ that are not regular points are \emph{singular points} of $K$. We say $K$ has smooth boundary if every point in the boundary $\partial K$ of $K\subseteq\Omega$ is a regular point of $K$. 

The following quantity is introduced by White \cite[Definition 4.2]{WhiTopolMCMCF}.

\begin{defn} \label{QDefn}
Let $\Omega$ be a closed subset of $\mathbb{R}^{n+1}$ with smooth boundary. For a closed set $K\subseteq\Omega$, we define $Q(K)$ to be the largest integer $l$ with the following properties: 
\begin{enumerate}
    \item The singular set $\mathrm{sing}(K)$ has Hausdorff dimension $\leq n-l$; \label{HausdorffItem}
    \item Let $p_i$ be a sequence of points in the interior of $K$ converging to a point $p$ in $\partial K$. Translate $K$ by $-p_i$ and dilate by $1/\mathrm{dist}(p_i,\partial K)$ to get $K_i$. Then a subsequence of the $K_i$ converges to a convex set $K^\prime$ (in $\mathbb{R}^{n+1}$ or in a closed half-space in $\mathbb{R}^{n+1}$) with smooth boundary, and the convergence is smooth on bounded sets;
    \item If $K^\prime$ is as in (2), then $\partial K^\prime$ has trivial $j$-th homotopy group for every $j<l$.
\end{enumerate}
If no such integer exists, let $Q(K)=-\infty$.
\end{defn}

It is shown in \cite[Proposition 4.3]{WhiTopolMCMCF} that for mean convex mean curvature flow of compact hypersurfaces, lower bounds on the quantity $Q$ are closely related to upper bounds on the Gaussian density. We adapt the arguments of \cite{WhiTopolMCMCF} to establish an analogous relationship between these bounds for strongly regular strictly monotone expander weak flows asymptotic to a cone.

\begin{lem} \label{LowBndQLem}
    Given $T\in \mathbb{R}$ and a $C^3$-regular cone $\mathcal{C}\subset\mathbb{R}^{n+1}$, let $\mathcal{S}=\{\Omega_t\}_{t\geq T}$ be a strongly regular strictly monotone expander weak flow asymptotic to $\mathcal{C}$ with starting time $T$. Assume that for $\mathcal{M}=\{\mathcal{H}^n\llcorner\partial\Omega_t\}_{t\geq T}$ and some $m\in [1,n-1]$
    \[
    \theta(\mathcal{M},X_0)\leq\lambda(\mathbb{S}^{m+1}\times\mathbb{R}^{n-m-1})
    \]
 at every point $X_0\in \mathbb{R}^{n+1}\times (T,\infty)$. Then there is an $R_1>1$ such that given $R>R_1$, if $K_t=\Omega_t\cap\bar{B}_R\subseteq\bar{B}_R$, then $Q(K_t)\geq m+1$ for each $t>T$.
\end{lem}

\begin{proof}
    Fix $\epsilon=10^{-10}$ and let $R_0$ be the radius given by Definition \ref{MonotoneExpanderFlowDefn}. We choose $R_1=4R_0$ and so, for $R>R_1$, the closed ball $\bar{B}_R$ has a collared neighborhood $U$ of $\partial B_R$ so that the restriction $\{K_t\}_{t\geq T}$ of $\mathcal{S}$ to $\bar{B}_R$ is regular in $U$. By our hypotheses and \eqref{EntSphereEqn}, every tangent flow at a given singularity is a self-similarly shrinking $\mathbb{S}^\ell\times\mathbb{R}^{n-\ell}$ with $\ell\geq m+1$. It follows from a dimension reduction argument (see \cite[\S 11]{WhiStrata}) that, for $t>T$, the singular set of $K_t$ has Hausdorff dimension at most $n-m-1$. 

    For $t>T$ fixed, let $p_i$ be a sequence of points in the interior of $K_t$ converging to a point $p\in \partial K_t$, and let  $\rho_i=1/\mathrm{dist}(p_i,\partial K_t)$ so $\rho_i\to\infty$. Translate $K_t$ by $-p_i$ and dilate by $\rho_i$ to get $K^i_t$. Then a subsequence of the $K^i_t$ converges to a set $K^\prime$. If $p\in\partial B_R$, then $p$ is a regular point of $K_t$. Thus, modulo a rigid motion, $K^\prime=\{x_1\geq 0\}\subset\mathbb{R}^{n+1}$ or $K^\prime=\{\mathbf{x}\cdot\mathbf{v}\geq 0,x_1\geq 0\}\subset\{x_1\geq 0\}$ for some $\mathbf{v}\in\mathbb{S}^n$, and the convergence is smooth on bounded sets. In particular, $K^\prime$ has trivial homotopy groups. So we may assume $p\in B_R$. Since the $K_t^i$ are part of a blow-up sequence of $\mathcal{S}$ and contain $B_1$, our hypotheses ensure $K^\prime$ is a convex subset of $\mathbb{R}^{n+1}$ with smooth boundary and the convergence is smooth on bounded sets. 
    
    It remains only to show that $\partial K^\prime$ has trivial $j$-th homotopy group for $j<m+1$. By convexity, either $\partial K^\prime$ is homeomorphic to $\mathbb{R}^n$, thus having trivial homotopy groups, or $\partial K^\prime$ is isometric to $K\times\mathbb{R}^{n-q}$ for some compact, convex set $K\subset\mathbb{R}^{q+1}$ -- see \cite[p. 3]{BusConvex}. As $\partial K^\prime$ is part of a limit flow at $P=(p,t)$, one appeals to a variant of Huisken's monotonicity formula (see \cite[\S 11]{WhiStrata}) and the Gaussian density bound to get
    \begin{equation} \label{EntLimitFlowEqn}
        \lambda(\partial K^\prime)\leq \theta(\mathcal{M},P)\leq \lambda(\mathbb{S}^{m+1}\times\mathbb{R}^{n-m-1})<2.
    \end{equation}
    In particular, this gives $q\geq 1$, as otherwise, $\partial K^\prime$ is the union of two parallel planes with entropy $2$. Thus, by results of Huisken \cite{HuiConvexMCF}, the mean curvature flow starting from $\partial K^\prime$ is given by the product of $\mathbb{R}^{n-q}$ and the flow $\mathcal{K}$ starting from  $\partial K\subset\mathbb{R}^{q+1}$, where $\mathcal{K}$ is smooth until it disappears in a round point. Thus, $\partial K^\prime$ is diffeomorphic to $\mathbb{S}^q\times\mathbb{R}^{n-q}$ and so, by Huisken monotonicity (see \cite{HuiAsympMCF}) and \eqref{EntLimitFlowEqn},
    \[
    \lambda(\mathbb{S}^q\times\mathbb{R}^{n-q}) \leq \lambda(\partial K^\prime) \leq\lambda(\mathbb{S}^{m+1}\times\mathbb{R}^{n-m-1}).
    \]
    Combined with \eqref{EntSphereEqn}, it follows that $q\geq m+1$ and so $\pi_j(\partial K^\prime)=\{0\}$ for $j<m+1$.
\end{proof}

Next we use results of White \cite{WhiTopolMCMCF} together with facts about algebraic topology to study topological properties of the expander flow as in Lemma \ref{LowBndQLem}.

\begin{prop} \label{ExpanderFlowTopolProp}
    Given $T\in \mathbb{R}$ and a $C^3$-regular cone $\mathcal{C}\subset\mathbb{R}^{n+1}$, let $\mathcal{S}=\{\Omega_t\}_{t\geq T}$ be a strongly regular strictly monotone expander weak flow asymptotic to $\mathcal{C}$ with starting time $T$. Assume that for $\mathcal{M}=\{\mathcal{H}^n\llcorner\partial\Omega_t\}_{t\geq T}$ and some $m\in [1,n-1]$
    \[
    \theta(\mathcal{M},X_0)\leq\lambda(\mathbb{S}^{m+1}\times\mathbb{R}^{n-m-1})
    \]
at every point $X_0\in \mathbb{R}^{n+1}\times (T,\infty)$. If $[a,b]\subset (T,\infty)$ and both $\Omega_a^c$ and $\Omega_b^c$ are path-connected, then the inclusion $\mathfrak{j}\colon \Omega_a^c\to\Omega_b^c$ induces homomorphisms 
    \[
    \mathfrak{j}_* \colon H_k(\Omega_a^c)\to H_k(\Omega_b^c)
    \]
    which are bijective when $k<m+1$ and surjective when $k=m+1$.  
        
    If, in addition, both $\partial\Omega_a$ and $\partial\Omega_b$ are smooth hypersurfaces, then the above remains true for the closures of $\Omega_a^c$ and $\Omega_b^c$, and the inclusion $\mathfrak{i}\colon\Omega_b\to\Omega_a$ induces homomorphisms
    \[
    \mathfrak{i}_* \colon H_k(\Omega_b)\to H_k(\Omega_a)
    \]
    which are bijective when $k>n-m-1$ and injective when $k=n-m-1$.
\end{prop}

\begin{proof}
     Fix $\epsilon=10^{-10}$ and let $R_0$ be the radius given by Definition \ref{MonotoneExpanderFlowDefn}. Let $R_1$ be the constant given by Lemma \ref{LowBndQLem}. For $R>4\max\{R_0,R_1\}$ fixed, we think of $\bar{B}_R$ as a manifold with boundary $S_R$. For $t\geq T$, let $K_t=\Omega_t\cap\bar{B}_R$. In what follows, denote by 
     $$
     K_t^c, \mathrm{int}(K_t), \mbox{ and }\partial K_t
     $$
the relative (in $\bar{B}_R$) complement, interior, and boundary, respectively, of $K_t$. Observe $K_t$ and $K_t^c$ are deformation retracts of $\Omega_t$ and $\Omega^c_t$, respectively. Thus, by homotopy equivalence, if the maps $\mathfrak{j}^\prime\colon K_a^c \to K_b^c$ and $\mathfrak{i}^\prime\colon K_b\to K_a$ are the inclusion maps, then it suffices to show:
    \begin{itemize}
    	\item[(i)] The induced homomorphisms $\mathfrak{j}^\prime_* \colon H_k(K_a^c)\to H_k(K_b^c)$ are bijective when $k<m+1$ and surjective when $k=m+1$; and;
    	\item [(ii)] Suppose both $\partial K_a$ and $\partial K_b$ are smooth, embedded manifolds with boundary and their boundaries meet $S_R$ transversally -- this last condition means $\partial K_a$ and $\partial K_b$ are manifolds with boundary in the sense of \cite[pg. 252]{HatAlgebraTopol}. Then the statement of (i) is still true for the relative (in $\bar{B}_R$) closures of $K_a^c$ and $K_b^c$, and the induced homomorphisms $\mathfrak{i}^\prime_* \colon H_k(K_b)\to H_k(K_a)$ are bijective when $k>n-m-1$ and injective when $k=n-m-1$.
    \end{itemize} 

    To see (i), observe that, by Lemma \ref{LowBndQLem}, $Q(K_t)\geq m+1$ for $t\in [a,b)$. Combined with our hypotheses, it follows from \cite[Theorem 5.2]{WhiTopolMCMCF} that $(K_b^c,K_a^c)$ is $(m+1)$-connected. Thus, by the relative Hurewicz theorem, e.g., \cite[Theorem 4.37]{HatAlgebraTopol}, as $K_b^c$ and $K_a^c$ are both path-connected, 
    \begin{equation} \label{VanishRelHomolEqn}
    H_\ell(K_b^c,K_a^c)=\{0\} \mbox{ for $\ell\leq m+1$}.
    \end{equation}
    The claim follows from the long exact sequence for $H_*(K_b^c,K_a^c)$ -- see \cite[pp. 115--117]{HatAlgebraTopol}.
    
    To see (ii), as $\partial K_a$ and $\partial K_b$ are smooth, embedded manifolds with boundary that meet $S_R$ transversally, it follows from homotopy equivalence that (i) holds for the relative (in $\bar{B}_R$) closures of $K_a^c$ and $K_b^c$. As arguing in \cite[Theorem 6.2]{WhiTopolMCMCF}, by the excision theorem (see \cite[Theorem 2.20]{HatAlgebraTopol}) and \eqref{VanishRelHomolEqn}, if $\ell\leq m+1$ then
    \[
    H_\ell(K_a\setminus \mathrm{int}(K_b),\partial K_a)\cong H_\ell(K_b^c,K_a^c)\cong \{0\}.
    \]
    Combined with the universal coefficients theorem (see \cite[Theorem 3.2]{HatAlgebraTopol}), one has
    \[
    H^\ell(K_a\setminus \mathrm{int}(K_b), \partial K_a)=\{0\}.
    \]
 The hypotheses on $K_a$ and $K_b$ ensure that $K_a\setminus \mathrm{int}(K_b)$ is an orientable compact manifold with boundary of dimension $n+1$ in the sense of \cite[pg. 252]{HatAlgebraTopol} -- we emphasize this is in a topological sense and not a smooth sense.  Moreover,  the boundary is $M\cup N$ where 
    \[
    M=\partial K_a \mbox{ and } N=\partial K_b \cup \left(\bigcup_{t\in [a,b)} \partial K_t\cap S_R\right),
    \]
   Thus, by Poincar\'{e}--Lefschetz duality  (see \cite[Theorem 3.43]{HatAlgebraTopol}), it follows that
    \[
    H_{n+1-\ell}(K_a\setminus \mathrm{int}(K_b),N)\cong H^\ell(K_a\setminus \mathrm{int}(K_b),M)\cong \{0\}.
    \]
    
    As $(K_a\setminus \mathrm{int}(K_b),N)$ is homotopy equivalent to $(K_a\setminus \mathrm{int}(K_b),\partial K_b)$, it follows that 
    \[
    H_{n+1-\ell}(K_a\setminus K_b,\partial K_b)=\{0\}.
    \]
    Hence, one may use the excision theorem and the long exact sequence for $H_*(K_a,K_b)$ as before to get the claim.
\end{proof}

We are now ready to prove Theorem \ref{NondegExpanderTopolThm}.

\begin{proof}[Proof of Theorem \ref{NondegExpanderTopolThm}]
    Suppose without loss of generality $\Gamma_-\neq\Gamma_+$, as otherwise, the inclusion of $\Omega_\pm(\Gamma_\pm)$ in $\Omega_\pm(\Gamma_\mp)$ is equal to the identity map and so induces isomorphisms on the homologies (that is, the theorem is trivially true). We first show there is a finite number of elements $\Gamma_1,\Gamma_2,\ldots,\Gamma_L$ of $\mathcal{E}(\mathcal{C})$ such that 
    \begin{itemize}
        \item $\Gamma_-=\Gamma_1\preceq\Gamma_2\preceq\cdots\preceq\Gamma_{L}=\Gamma_+$;
        \item The $\Gamma_\ell$ alternate between being (strictly) stable and unstable;
        \item If $\Gamma_\ell$ is unstable, then there exist two perturbations by normal graphs of $\Gamma_\ell$ to its either side, $\Gamma_\ell^-\preceq\Gamma_\ell\preceq\Gamma_\ell^+$ such that $\Omega_\pm(\Gamma_\ell^\pm)$ is strictly expander mean convex, and there is a strongly regular strictly monotone expander weak flow $\mathcal{S}_\ell$ asymptotic to $\mathcal{C}$ that starts from $\Omega_\pm(\Gamma_\ell^\pm)$ and converges, in $C^\infty_{loc}(\mathbb{R}^{n+1})$, to $\Omega_\pm(\Gamma_{\ell\pm1})$ as $t\to\infty$ provided $1\leq \ell\pm 1\leq L$.
    \end{itemize}
    
    To see this, first observe that every element of $\mathcal{E}(\mathcal{C})$ is connected because $\mathcal{L}(\mathcal{C})$ is connected and there are no closed expanders. By \cite[Proposition 6.3]{BWTopolUnique} and \cite[Proposition 3.2]{BWMinMax}, if $\Gamma_\pm$ is unstable, then there is a $C^3$-asymptotically conical hypersurface $\Gamma'_\pm$ given by a normal graph of $\Gamma_\pm$ that is trapped between $\Gamma_-$ and $\Gamma_+$, and so that $\Omega_\mp(\Gamma'_\pm)$ is strictly expander mean convex. Thus, by Proposition \ref{MonotoneExpanderFlowProp}, there is a strongly regular strictly monotone expander weak flow asymptotic to $\mathcal{C}$ that starts from $\Omega_\mp(\Gamma'_\pm)$ and, as $t\to\infty$, converges smoothly to $\Omega_\mp(\hat{\Gamma}_\pm)$ for stable $\hat{\Gamma}_\pm\in\mathcal{E}(\mathcal{C})$. And the construction ensures $\Gamma_-\preceq\hat{\Gamma}_\pm\preceq\Gamma_+$. Thus, we may assume both $\Gamma_+$ and $\Gamma_-$ are stable. In what follows, we prove the claim by extending the arguments of \cite{BWTopolUnique} for $m=1$ to general $m$. 
    
    Let $\mathcal{E}_S(\Gamma_-,\Gamma_+)$ be the set of elements $\Gamma$ of $\mathcal{E}(\mathcal{C})$ that are stable and such that $\Gamma_-\preceq\Gamma\preceq\Gamma_+$. By compactness properties of the space of stable expanders -- i.e., \cite[Proposition 4.4]{BWTopolUnique} -- and the nondegeneracy hypothesis on $\cC$, $\mathcal{E}_S(\Gamma_-,\Gamma_+)$ is a finite set. We use induction on the number of elements $J=|\mathcal{E}_S(\Gamma_-,\Gamma_+)|$.  There are at least two elements by hypothesis and, when $J=2$, as $\Gamma_+,\Gamma_-\in\mathcal{E}(\mathcal{C})$ are (strictly) stable and $\Gamma_-\preceq\Gamma_+$, a min-max construction (see \cite{BWMinMax}) yields a $\tilde{\Gamma}\in\mathcal{E}(\mathcal{C})\setminus\{\Gamma_+,\Gamma_-\}$ with $\Gamma_-\preceq\tilde{\Gamma}\preceq\Gamma_+$. Indeed, our hypotheses ensure that $\tilde{\Gamma}$ is connected and unstable, and so a straightforward modification of the arguments in the previous paragraph gives the claim for $J=2$. Next, suppose $\tilde{\Gamma}_+$ is a maximal element of $\mathcal{E}_S(\Gamma_-,\Gamma_+)\setminus\{\Gamma_+\}$, that is, there are no $\Gamma\in\mathcal{E}_S(\Gamma_-,\Gamma_+)\setminus\{\Gamma_+,\tilde{\Gamma}_+\}$ such that $\tilde{\Gamma}_+\preceq\Gamma\preceq\Gamma_+$; such an element exists as we are considering a finite ordered set. For $J\geq 3$, one has $|\mathcal{E}_S(\tilde{\Gamma}_+,\Gamma_+)|=2$ and $|\mathcal{E}_S(\Gamma_-,\tilde{\Gamma}_+)|\leq J-1$. The claim follows from the induction hypothesis.
    
    For $1\leq \ell\leq L-1$, let $\mathfrak{i}^+_\ell$ be the inclusion of $\Omega_+(\Gamma_{\ell+1})$ in $\Omega_+(\Gamma_\ell)$ and $\mathfrak{i}^-_\ell$ the inclusion of $\Omega_-(\Gamma_{\ell})$ in $\Omega_-(\Gamma_{\ell+1})$. If $\Gamma_\ell$ is unstable, one can ensure, by \cite[Lemma 3.5]{BWProper},
    $$
    \lambda(\Gamma_\ell^\pm)<\lambda(\Gamma_\ell)+\epsilon=\lambda(\mathcal{C})+\epsilon<\lambda(\mathbb{S}^m\times\mathbb{R}^{n-m}).
    $$
 By Lemma \ref{GaussDenLem}, if $\mathcal{M}_\ell$ is the measure flow associated with $\mathcal{S}_\ell$ -- see Item (5) of Definition \ref{MonotoneExpanderFlowDefn} -- 
 one has $\theta(\mathcal{M}_\ell,X_0)\leq\lambda(\mathbb{S}^{m+1}\times\mathbb{R}^{n-m-1})$ at every spacetime point $X_0$. Thus, Proposition \ref{ExpanderFlowTopolProp} together with homotopy equivalence implies that the induced homomorphisms $(\mathfrak{i}^+_\ell)_* \colon H_k(\Omega_+(\Gamma_{\ell+1}))\to H_k(\Omega_+(\Gamma_{\ell}))$ are bijective when $n-m\leq k\leq m$, are injective when $n-m-1=k\leq m$, and are surjective when $n-m\leq k=m+1$. The same is true for the induced maps $(\mathfrak{i}^-_\ell)_*\colon H_k(\Omega_-(\Gamma_\ell))\to H_k(\Omega_-(\Gamma_{\ell+1}))$. As $\mathfrak{i}^+=\mathfrak{i}^+_1\circ\cdots\circ\mathfrak{i}^+_{L-1}$ and $\mathfrak{i}^-=\mathfrak{i}^-_{L-1}\circ\cdots\circ\mathfrak{i}^-_1$, the claim follows.
\end{proof}

Next we aim to extend Theorem \ref{NondegExpanderTopolThm} by dropping the nondegeneracy hypothesis. Because nondegenerate cones are generic in an appropriate space of cones, one would like to perturb the pair $\Gamma_+,\Gamma_-\in\mathcal{E}(\mathcal{C})$ with $\Gamma_-\preceq\Gamma_+$ and cone $\cC$ in order to produce expander mean-convex hypersurfaces $\Gamma_+^\prime,\Gamma_-^\prime$ asymptotic to $\cC'$ where $\mathcal{C}^\prime$ is nondegenerate -- see \cite[Corollary 1.2]{BWBanach}. Indeed, Propositions \ref{MonotoneExpanderFlowProp} and \ref{ExpanderFlowTopolProp} together yield $\hat{\Gamma}_+,\hat{\Gamma}_-\in\mathcal{E}(\mathcal{C}^\prime)$ along with relationships between the  topologies of $\Gamma_\pm^\prime$ and $\hat{\Gamma}_\pm$.  However, to complete this argument one must be able to control the relationship between the partial orders and direction of the expander mean curvature vectors which is a subtle issue. Specifically, one must also have $\hat{\Gamma}_+\preceq\Gamma_+^\prime$, $\Gamma_-^\prime\preceq\hat{\Gamma}_-$, and $\hat{\Gamma}_-\preceq\hat{\Gamma}_+$, which do not, in general, hold.  To overcome this, we prove a restricted result by requiring $\Gamma_+$ and $\Gamma_-$ to be the greatest and least element of $\mathcal{E}(\mathcal{C})$ -- whose existence is shown in \cite[Theorem 4.1]{BWTopolUnique}.  

\begin{prop} \label{DegExpanderTopolProp}
    For $2\leq n\leq 6$, let $\mathcal{C}$ be a degenerate $C^5$-regular cone in $\mathbb{R}^{n+1}$ such that $\mathcal{L}(\mathcal{C})$ is connected and, for some $m\in [1,n-1]$,
    \[
    \lambda(\mathcal{C})<\lambda(\mathbb{S}^{m}\times\mathbb{R}^{n-m}).
    \]
    If $\Gamma_+$ and $\Gamma_-$ are respectively the greatest and least element of $\mathcal{E}(\mathcal{C})$, i.e., $\Gamma_-\preceq\Gamma\preceq\Gamma_+$ for $\Gamma\in\mathcal{E}(\mathcal{C})$, then the inclusions $\mathfrak{i}^\pm \colon \Omega_\pm(\Gamma_\pm) \to \Omega_\pm(\Gamma_\mp)$ induce homomorphisms
    \[
    \mathfrak{i}^\pm_* \colon H_k(\Omega_\pm(\Gamma_\pm)) \to H_k(\Omega_\pm(\Gamma_\mp))
    \]
    such that 
    \begin{enumerate}
        \item When $n-m\leq k\leq m$, the maps $\mathfrak{i}^\pm_*$ are bijective;
        \item When $n-m-1=k\leq m$, the maps $\mathfrak{i}^\pm_*$ are injective;
        \item When $n-m\leq k=m+1$, the maps $\mathfrak{i}^\pm_*$ are surjective.
    \end{enumerate}
\end{prop}

\begin{proof}
    First we show there is a sequence of nondegenerate $C^4$-regular cones $\mathcal{C}_\ell\subset\mathbb{R}^{n+1}$ with link $\mathcal{L}(\mathcal{C}_\ell)$ so that $\mathcal{L}(\mathcal{C}_\ell)\to\mathcal{L}(\mathcal{C})$ in $C^4(\mathbb{S}^n)$ and $\Gamma^\ell_+\to\Gamma_+$ in $C^\infty_{loc}(\mathbb{R}^{n+1})$, where $\Gamma_+^\ell$ is the greatest element of $\mathcal{E}(\mathcal{C}_\ell)$. In particular, for sufficiently large $\ell$,
    \begin{equation} \label{EntConeEqn}
    \lambda(\mathcal{C}_\ell)<\lambda(\mathcal{C})+\epsilon<\lambda(\mathbb{S}^m\times\mathbb{R}^{n-m}).
    \end{equation}
    To see this, first observe by \cite[Theorem 4.1]{BWTopolUnique} one has $\Gamma_+$ is a stable self-expander. Let  $\mathcal{K}$ denote the space of Jacobi fields on $\Gamma_+$ that fix the infinity.  Using standard spectral theory and the results of \cite{BerAsympEigenfunc} -- see also \cite[Lemma 6.1]{BWBanach} and \cite[Proposition 3.2]{BWMinMax} -- together with the connectedness of  $\Gamma_+$ it follows that $\dim \mathcal{K}\leq 1$ and every element of $\mathcal{K}$ either is identically $\mathbf{0}$ or has a sign and an asymptotic expansion with leading order term of the form $\alpha r^{-n-1} e^{-r^2/4}$, where $r=|\mathbf{x}|$ and $\alpha$ is a function on $\mathcal{L}(\mathcal{C})$. 
       
    We now appeal to results of \cite{BWBanach}  to perturb the asymptotic cone and the self-expander $\Gamma_+$.  First, recall that $\mathcal{ACH}^4_n(\Gamma_+)$ denotes the space of $C^4$-asymptotically conical embeddings of $\Gamma_+$ into $\mathbb{R}^{n+1}$, equipped with the weighted $C^4$ norm with weight $r$ -- see \cite[\S 3.2]{BWBanach} for the precise definition. Following \cite{BWBanach}, say a $C^3$-regular cone $\cC'$ is nondegenerate if no self-expander asymptotic to $\cC'$ admits a non-trivial Jacobi field -- as shown in \cite{BWBanach} such cones are generic in a suitable sense. With this in mind we use \cite[Theorem 7.1]{BWBanach} to choose a sequence of $C^4$-regular cones $\mathcal{C}_\ell$ along with $C^4$-asymptotically conical hypersurfaces $\Gamma_\ell$ with nondegenerate asymptotic cone so that the $\cC_{\ell}$ converge to $\cC$ and the $\Gamma_\ell$ converge  $\Gamma_+$ in $\mathcal{ACH}^4_n(\Gamma_+)$. Here we understand the convergence by identifying the $\Gamma_\ell$ with the parameterizations given by the inverse of the nearest point projections onto $\Gamma_+$ which are elements of $\mathcal{ACH}^4_n(\Gamma_+)$.
    
    It follows from \cite[Theorem 7.1]{BWBanach} that the cones $\cC_{\ell}$ and hypersurfaces $\Gamma_\ell$ can be chosen so that either:
    \begin{enumerate}[(i)]
    	\item \label{Case1} If $\mathcal{K}=\{\mathbf{0}\}$, then $\Gamma_\ell\in\mathcal{E}(\mathcal{C}_\ell)$ for each $\ell$; or;
    	\item \label{Case2} If $\mathcal{K}\neq \{\mathbf{0}\}$, then $2\mathbf{H}_{\Gamma_\ell}-\mathbf{x}^\perp$ is nowhere vanishing and the pushforward of it by the nearest projection onto $\Gamma_+$ is a nonzero element of $\mathcal{K}$.
    \end{enumerate}
If case \eqref{Case2} occurs, then one can modify the choice of cones in order to ensure that $2\mathbf{H}_{\Gamma_\ell}-\mathbf{x}^\perp$ points into $\Omega_+(\Gamma_\ell)$ -- here one chooses the boundary link of $\mathcal{L}(\mathcal{C}_\ell)$ in the obvious way compatible with that of $\mathcal{L}(\cC)$.  This ensures that for the choice of unit normal $\mathbf{n}_{\Gamma_\ell}$ on $\Gamma_{\ell}$ that is compatible with the boundary link -- i.e., that points out of $\Omega_+(\Gamma_\ell)$ -- one has $2H_{\Gamma_\ell}+\mathbf{x}\cdot\mathbf{n}_{\Gamma_\ell}>0$. 
    
    In case \eqref{Case1}, it is readily seen that the construction and the compactness properties of spaces of expanders \cite[Proposition 4.4]{BWTopolUnique} ensure that the $\Gamma^\ell_+$ converge smoothly on compact sets to $\Gamma_+$. In case \eqref{Case2}, by Proposition \ref{MonotoneExpanderFlowProp}, the expander flow starting from $\Gamma_\ell$ is contained in $\Omega_+(\Gamma_\ell)$ and is asymptotic at time $\infty$ to a stable self-expander $\hat{\Gamma}_\ell\in\mathcal{E}(\mathcal{C}_\ell)$. As $\Gamma_\ell\preceq\hat{\Gamma}_\ell$ and $\Gamma_\ell\to\Gamma_+$, invoking compactness again gives that the $\hat{\Gamma}_\ell$ converge, in $C^\infty_{loc}(\mathbb{R}^{n+1})$, to a limit expander $\hat{\Gamma}\in \mathcal{E}(\cC)$ that satisfies $\Gamma_+\preceq \hat{\Gamma}$.  As $\Gamma_+$ is the greatest element one must have $\Gamma_+=\hat{\Gamma}$.  The same reasoning ensures that, up to passing to a further subsequence,  the greatest elements, $\Gamma_+^\ell$, of $\mathcal{E}(\cC_\ell)$ also converge to $\Gamma_+$. This proves the claim.
    
    Using the sequence of cones, $\cC_\ell$, from the previous step and making the same argument with $\Gamma_-$ in place of $\Gamma_+$ yields a sequence of $C^4$-asymptotically conical hypersurfaces $\Upsilon_\ell$ with asymptotic cone $\mathcal{C}_\ell$ that converges to $\Gamma_-$ in $\mathcal{ACH}^4_n(\Gamma_-)$ and such that, up to passing to a subsequence, one of the following holds: 
    \begin{enumerate}
        \item[(a)] $\Upsilon_\ell\in\mathcal{E}(\mathcal{C}_\ell)$ for each $\ell$;
        \item[(b)] $2H_{\Upsilon_\ell}+\mathbf{x}\cdot\mathbf{n}_{\Upsilon_\ell}<0$ for each $\ell$;
        \item[(c)] $2H_{\Upsilon_\ell}+\mathbf{x}\cdot\mathbf{n}_{\Upsilon_\ell}>0$ for each $\ell$.
    \end{enumerate}
    Here $\mathbf{n}_{\Upsilon_\ell}$ is the outward normal to $\Omega_+(\Upsilon_\ell)$.  Observe that three cases may occur rather than two as in the previous step as we have already fixed the sequence of cones.
    
    When case (a) or (b) occurs, one argues as before to see the least elements $\Gamma^\ell_-$ of $\mathcal{E}(\mathcal{C}_\ell)$ converge smoothly on compact sets to $\Gamma_-$. By \cite[Proposition 3.3]{BWProper} there is a radius $R>1$ so that $\Omega_\pm(\Gamma^\ell_-)\cap\bar{B}_R$ is a deformation retract of $\Omega_\pm(\Gamma^\ell_-)$ for each $\ell$. The nature of convergence ensures that $\Omega_\pm(\Gamma^\ell_-)$ is homotopy equivalent to $\Omega_\pm(\Gamma_-)$ for sufficiently large $\ell$. Likewise, $\Omega_\pm(\Gamma^\ell_+)$ is homotopy equivalent to $\Omega_\pm(\Gamma_+)$ for large $\ell$. The claim follows from this and Theorem \ref{NondegExpanderTopolThm} for $\Gamma_-^\ell\preceq\Gamma^\ell_+$ (in view of \eqref{EntConeEqn}).
    \begin{figure}
    	 \centering
    	\resizebox{4in}{!}{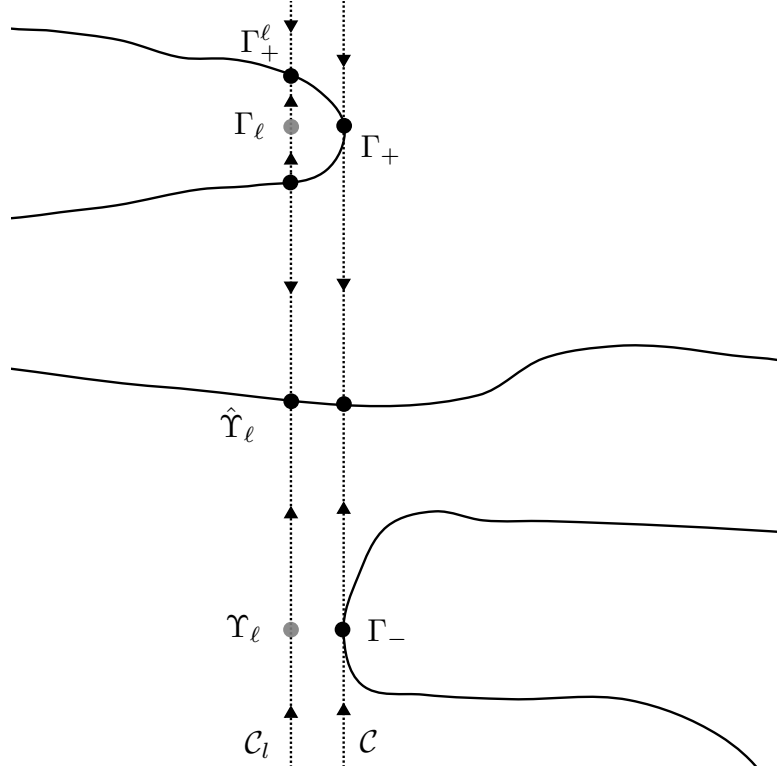}
    	\caption{A schematic illustration of the situation in case \eqref{Case2} together with case (c).  The horizontal axis is space of cones, the vertical axis the space of hypersurfaces asymptotic to the given cone where height corresponds to the order, $\preceq$. The arrows represent directions of flow lines.}
    \label{NonGenCone2c}
    \end{figure}
    
    It remains only to deal with case (c).  We refer to Figure \ref{NonGenCone2c}. Fix a large $\ell$ so that $\Omega_\pm(\Gamma_+^\ell)$ and $\Omega_\pm(\Upsilon_\ell)$ are homotopy equivalent to $\Omega_\pm(\Gamma_+)$ and $\Omega_\pm(\Gamma_-)$, respectively. Observing $\Upsilon_\ell\preceq\Gamma_+^\ell$, we then let $\mathfrak{i}_\ell^+\colon \Omega_+(\Gamma_+^\ell)\to\Omega_+(\Upsilon_\ell)$ and $\mathfrak{i}_\ell^-\colon \Omega_-(\Upsilon_\ell)\to\Omega_-(\Gamma_+^\ell)$ be the inclusion maps. Hence, it suffices to show the induced maps $(\mathfrak{i}_\ell^+)_* \colon H_k(\Omega_+(\Gamma_+^\ell))\to H_k(\Omega_+(\Upsilon_\ell))$ and $(\mathfrak{i}^-_\ell)_* \colon H_k(\Omega_-(\Upsilon_\ell))\to H_k(\Omega_-(\Gamma^\ell_+))$ are bijective when $n-m\leq k\leq m$, are injective when $n-m-1=k\leq m$, and are surjective when $n-m\leq k=m+1$. 
    
    By Proposition \ref{MonotoneExpanderFlowProp}, there is a strongly regular strictly monotone expander weak flow $\mathcal{S}_\ell$ asymptotic to $\mathcal{C}_\ell$ with initial data $\Omega_+(\Upsilon_\ell)$, and the flow converges smoothly to a closed set with smooth boundary $\hat{\Upsilon}_\ell\in\mathcal{E}(\mathcal{C}_\ell)$ satisfying ${\Upsilon}_{\ell}\preceq \hat{\Upsilon}_\ell\preceq \Gamma_+^\ell$. In particular, for large $t$, the time-$t$ slice of the flow line is smoothly isotopic to $\hat{\Upsilon}_\ell$. By construction 
    \[
    \lambda(\Upsilon_\ell)<\lambda(\Gamma_-)+\epsilon=\lambda(\mathcal{C})+\epsilon<\lambda(\mathbb{S}^m\times\mathbb{R}^{n-m})
    \]
    and so the Huisken monotonicity \cite{HuiAsympMCF} implies, for $X_0\in \mathbb{R}^{n+1}\times (0,\infty)$, that
    \[
    \theta(\mathcal{M}_\ell,X_0)\leq\lambda(\mathbb{S}^{m+1}\times\mathbb{R}^{n-m-1})
    \]
where $\mathcal{M}_\ell$ is the measure flow associated with $\mathcal{S}_\ell$ -- see Item (5) of Definition \ref{MonotoneExpanderFlowDefn}. The result follows by combining Proposition \ref{ExpanderFlowTopolProp} applied to $\mathcal{S}_\ell$ which yields relations between the topologies of $\Upsilon_\ell\preceq \hat{\Upsilon}_\ell$ together with Theorem \ref{NondegExpanderTopolThm}  which gives relationships between the topologies of  $\hat{\Upsilon}_\ell\preceq\Gamma^\ell_+$.  Observe that Theorem \ref{NondegExpanderTopolThm} applies to $\hat{\Upsilon}_\ell\preceq \Gamma_\ell^+$ because of \eqref{EntConeEqn} and because the $\cC_\ell$ were chosen to be nondegenerate.
\end{proof}

\section{Density of minimal cones} \label{DenMinimalConeSec}
In this section we prove the main results about densities of minimal cones. We continue to use the conventions of Section \ref{TopolChangeSec}. We will repeatedly invoke the following auxiliary lemma about properties of the greatest and least element of $\mathcal{E}(\mathcal{C})$ for a minimal cone $\mathcal{C}$ -- see \cite[Theorem 4.1]{BWTopolUnique} for the existence of those elements.

\begin{lem} \label{StarShapeExpanderLem}
    For $3\leq n\leq 6$, let $\mathcal{C}\subset\mathbb{R}^{n+1}$ be a non-flat regular minimal cone. If $\Sigma_+$ and $\Sigma_-$ are respectively the greatest and least element of $\mathcal{E}(\mathcal{C})$, then the following is true: 
    \begin{enumerate}
    	\item $\mathcal{L}(\cC)$ is connected;
        \item $\Sigma_\pm\subset\mathrm{int}(\Omega_\pm(\mathcal{C}))$; 
        \item The projections $\Pi^\pm\colon\Sigma_\pm\to\mathbb{S}^n$ given by
        \[
        \Pi^\pm(p)=\frac{\mathbf{x}(p)}{|\mathbf{x}(p)|} 
        \]
        are diffeomorphisms onto their images  $\mathrm{int}(\Omega_\pm(\mathcal{C}))\cap\mathbb{S}^n$; 
        \item The sets $\Omega_\mp(\Sigma_\pm)$ are star-shaped relative to $\mathbf{0}$;
        \item The inclusions of $\Sigma_\pm$ into $\mathrm{int}(\Omega_\pm(\cC))$ are deformation retracts.
    \end{enumerate}
\end{lem}

\begin{proof}
	First of all, because $\mathcal{L}(\cC)$ is a smooth minimal hypersurface in $\mathbb{S}^n$, the Frankel theorem \cite{Frankel} ensures it has only one connected component.
	
    As $\mathcal{C}$ is non-flat regular minimal cone, it is $E$-stationary and singular.   It follows from the regularity theory for area-minimizing hypersurfaces that $\cC$ cannot be $E$-minimizing. Thus, $\Sigma_\pm\subset\mathrm{int}(\Omega_\pm(\mathcal{C}))$, as otherwise, one may use a minimization procedure (see \cite[pp. 13--14]{IlmNoteMCF} and \cite[\S 6]{DinMCExpander} or \cite[\S 4]{BWTopolUnique}) to construct a $\Sigma^\prime_\pm\in\mathcal{E}(\mathcal{C})$ so that $\Sigma^\prime_\pm$ is contained in the interior of $\Omega_\pm(\Sigma_\pm)\cap\Omega_\pm(\mathcal{C})$,  either contradicting $\Sigma_+$ being the greatest element or $\Sigma_-$ being the least.
    
    Next we show $\sqrt{t_1}\Sigma_+\preceq\sqrt{t_2}\Sigma_+$ for $0<t_1<t_2$, which implies $\mathbf{H}_{\Sigma_+}$ vanishes nowhere and points into $\Omega_+(\Sigma_+)$. To see this, we let $\delta=t_2-t_1$ and consider two mean curvature flows $\mathcal{M}=\{\sqrt{t}\Sigma_+\}_{t>0}$ and $\mathcal{M}^\prime=\{\sqrt{t+\delta}\Sigma_+\}_{t>0}$. As $\Sigma_+$ is $C^2$-asymptotic to $\mathcal{C}$, there is a radius $R>1$ so that $\sqrt{t+\delta}\Sigma_+\setminus B_{R}$ can be written as an almost flat normal graph over (a subset of) $\sqrt{t}\Sigma_+$ for $0<t\leq t_1$. As $\Omega_+(\sqrt{\delta}\Sigma_+)\subset\mathrm{int}(\Omega_+(\mathcal{C}))$ the continuity of flows ensures that, for $t>0$ small,  $\Omega_+(\sqrt{t+\delta}\Sigma_+)\cap \partial B_{2R}$ is contained in the interior of $\Omega_+(\sqrt{t}\Sigma_+)\cap\partial B_{2R}$. Suppose $s$ is the first time in $(0,t_1]$ such that this fails. Then the strict maximum principle on noncompact regions implies the restriction of $\mathcal{M}$ to $(\mathbb{R}^{n+1}\setminus \bar{B}_{2R})\times (0,s]$ is disjoint from $\mathcal{M}^\prime$. Thus, there is a time $s^\prime\in (0,s]$ so that the restriction of $\mathcal{M}$ to $(0,s^\prime)$ lies on one side of $\mathcal{M}^\prime$ and $\mathcal{M}$ touches $\mathcal{M}^\prime$ at a point in $\bar{B}_{2R}$ at time $s^\prime$, which violates the strict maximum principle on compact regions. That is, $\sqrt{t+\delta}\Sigma_+\cap\sqrt{t}\Sigma_+\cap\partial B_{2R}=\emptyset$ for $t\in (0,t_1]$. The claim follows from the strict maximum principle.
    
    Swapping the orientation, the arguments above imply $\mathbf{H}_{\Sigma_-}$ vanishes nowhere and points into $\Omega_-(\Sigma_-)$. The expander equation \eqref{ExpanderEqn} implies $\mathbf{x}\cdot\mathbf{n}_{\Sigma_+}>0$ and $\mathbf{x}\cdot\mathbf{n}_{\Sigma_-}<0$, where $\mathbf{n}_{\Sigma_\pm}$ is the outward normal to $\Omega_+(\Sigma_\pm)$. It follows from \cite[Proposition 5.1]{BWACShrinker} that the projection $\Pi^\pm\colon \Sigma_\pm\to\mathbb{S}^n$ is a diffeomorphism onto its image, and the image is $\mathrm{int}(\Omega_\pm(\mathcal{C}))\cap\mathbb{S}^n$. This also implies that $\Omega_\mp(\Sigma_\pm)$ is star-shaped relative to $\mathbf{0}$ and that the inclusions of $\Sigma_\pm$ into $\mathrm{int}(\Omega_\pm(\cC))$ are deformation retracts and so completes the proof.
\end{proof}

Theorem \ref{NonContractThm} follows from the proposition below.

\begin{prop} \label{ContractProp}
    For $3\leq n\leq 6$, let $\mathcal{C}$ be a regular minimal cone in $\mathbb{R}^{n+1}$ with 
    \[
    \lambda(\mathcal{C})<\lambda(\mathbb{S}^{n-1}\times\mathbb{R}).
    \]
    Then, for $\Sigma\in\mathcal{E}(\mathcal{C})$, both $\Omega_+(\Sigma)$ and $\Omega_-(\Sigma)$ are contractible. This means $\Sigma$ and both components of $\mathbb{S}^n\setminus\mathcal{L}(\mathcal{C})$ are contractible.  
\end{prop}
\begin{rem}
	The fact that the two components of $\mathbb{S}^n\setminus \mathcal{L}(\cC)$ are contractible implies $\mathcal{L}(\cC)$ is a homology sphere, but can not rule out the presence of torsion elements in the fundamental group when $n\geq 4$.
\end{rem}

\begin{proof}[Proof of Proposition \ref{ContractProp}]
    By the maximum principle, if $\mathcal{C}$ is flat, then the only self-expander asymptotic to $\mathcal{C}$ is a hyperplane and the claim holds trivially. Otherwise, let $\Sigma_+$ and $\Sigma_-$ be respectively the greatest and least element of $\mathcal{E}(\mathcal{C})$. By Lemma \ref{StarShapeExpanderLem}, $\mathcal{L}(\cC)$ is is connected and both $\Omega_-(\Sigma_+)$ and $\Omega_+(\Sigma_-)$ are star-shaped relative to $\mathbf{0}$ and so are contractible. Combined with \cite[Theorem 1.1]{BWTopolUnique}, it follows that, for $\Sigma\in\mathcal{E}(\mathcal{C})$, both $\Omega_-(\Sigma)$ and $\Omega_+(\Sigma)$ are contractible. Furthermore, by Lemma \ref{StarShapeExpanderLem}, if $\omega_\pm=\mathrm{int}(\Omega_\pm(\mathcal{C}))\cap\mathbb{S}^n$, then both $\Sigma_\pm$ and $\omega_\pm$ are homotopy equivalent to $\Omega_\pm(\Sigma_\pm)$, and thus both $\Sigma_\pm$ and $\omega_\pm$ are contractible. Finally, by \cite[Theorem 1.1]{BWTopolUnique} every $\Sigma\in\mathcal{E}(\mathcal{C})$ is diffeomorphic to $\Sigma_+$ and so is contractible.
\end{proof}

Theorem \ref{NonHomolBallThm} follows from the proposition below.

\begin{prop} \label{HomolBallProp}
    For $3\leq n\leq 6$, let $\mathcal{C}$ be a regular minimal cone in $\mathbb{R}^{n+1}$ with
    \[
    \lambda(\mathcal{C})<\lambda(\mathbb{S}^{n-2}\times\mathbb{R}^2).
    \]
    If $\Sigma_+$ and $\Sigma_-$ are respectively the greatest and least element of $\mathcal{E}(\mathcal{C})$, then, for $k\geq 0$,
    \[
    \tilde{H}_k(\Omega_\mp(\Sigma_\pm))=\tilde{H}_k(\Omega_\pm(\Sigma_\pm))=\{0\}.
    \]
    This means $\Sigma_+$, $\Sigma_-$, and the two components of $\mathbb{S}^n\setminus\mathcal{L}(\mathcal{C})$ are all homology balls. In particular, $\mathcal{L}(\mathcal{C})$ is a homology sphere and, when $n=3$, the cone $\mathcal{C}$ is flat.
\end{prop}

\begin{proof}
   If $\mathcal{C}$ is flat, then so is $\Sigma_\pm$. Otherwise, by Lemma \ref{StarShapeExpanderLem}, $\mathcal{L}(\cC)$ is connected and $\Omega_\mp(\Sigma_\pm)$ are star-shaped relative to $\mathbf{0}$. In particular, it follows that $\Omega_\mp(\Sigma_\pm)$ are homotopy equivalent to $B^{n+1}_1$ and so have trivial reduced homology groups.
    
    Likewise, if $\omega_\pm=\mathrm{int}(\Omega_\pm(\mathcal{C}))\cap\mathbb{S}^n$, then $\Omega_\pm(\Sigma_\pm)$ is homotopy equivalent to $\omega_\pm$ and so, for $k\geq 0$,
    \[
        H_k(\Omega_\pm(\Sigma_\pm))\cong H_k(\omega_\pm).
    \]
    In particular, by the Alexander duality theorem -- e.g.,  \cite[Theorem 3.44]{HatAlgebraTopol} -- one immediately has
    \[
    H_{n+1}(\Omega_\pm(\Sigma_\pm))\cong H_n(\Omega_\pm(\Sigma_\pm)) \cong H_{n-1}(\Omega_\pm(\Sigma_\pm))\cong \{0\}.
    \]
    
    By Theorem \ref{NondegExpanderTopolThm} and Proposition \ref{DegExpanderTopolProp}, for $2\leq k\leq n-2$, one has $H_k(\Omega_\pm(\Sigma_\pm))$ is isomorphic to $H_k(\Omega_\pm(\Sigma_\mp))$ and so is trivial (this is vacuous for $n=3$). Moreover, for $1=k\leq n-2$, one has $H_1(\Omega_\pm(\Sigma_\pm))$ is isomorphic to a subgroup of $H_1(\Omega_\pm(\Sigma_\mp))$. As $H_1(\Omega_\pm(\Sigma_\mp))$ is trivial so is $H_1(\Omega_\pm(\Sigma_\pm))$. It follows that $\Omega_\pm(\Sigma_\pm)$ is a homology ball. Furthermore, by Lemma \ref{StarShapeExpanderLem} both $\Sigma_\pm$ and $\omega_\pm$ are homotopy equivalent to $\Omega_\pm(\Sigma_\pm)$ and so are homology balls.
	
	To complete the proof, observe that, by the Mayer--Vietoris sequence, e.g.,  \cite[p. 149]{HatAlgebraTopol}, as both $\omega_+$ and $\omega_-$ are homology balls, $\mathcal{L}(\mathcal{C})$ is a homology sphere. When $n=3$, the link $\mathcal{L}(\cC)$ is a surface that is a homology sphere and so, by the classification of surfaces, must be a topological $\mathbb{S}^2$. That is, $\mathcal{L}(\cC)$ is a minimal two-sphere in $\mathbb{S}^3$ and so, by results of Almgren \cite{AlmMinSurfSphere}, must be totally geodesic; in this case $\mathcal{C}$ is flat.
\end{proof}
We also obtain a weaker topological restriction under weaker entropy bounds.
\begin{prop}\label{WeakProp}
	For $5\leq n\leq 6$, let $\mathcal{C}$ be a regular minimal cone in $\mathbb{R}^{n+1}$ with
	\[
	\lambda(\mathcal{C})<\lambda (\mathbb{S}^{n-3}\times \Real^{3}).
	\]
	If $\Sigma_+$ and $\Sigma_-$ are respectively the greatest and least elements of $\mathcal{E}(\mathcal{C})$, then 	the following is true:
	\begin{enumerate}
	    \item $\Omega_\mp(\Sigma_\pm)$ is a homology ball;
	    \item When $k\neq 1, n-2$, one has $\tilde{H}_k(\Omega_\pm(\Sigma_\pm))=\{0\}$. Otherwise, $\tilde{H}_k(\Omega_\pm(\Sigma_\pm))$ is a torsion-free group.
	\end{enumerate}
    This means that when $k\neq 1,n-2$, one has $\tilde{H}_k(\Sigma_\pm)=\{0\}$ and otherwise, $\tilde{H}_k(\Sigma_\pm)$ is a torsion-free group. The same is true for the two components of $\mathbb{S}^n\setminus \mathcal{L}(\mathcal{C})$. In particular,  $H_k(\mathcal{L}(\cC))=\{0\}$ when  $2\leq k\leq n-3$ and both $H_1(\mathcal{L}(\cC))$ and $H_{n-2}(\mathcal{L}(\mathcal{C}))$ are torsion-free groups. 
\end{prop}

\begin{rem}
Observe that this result applies to $\cC_{2,2}\subset \mathbb{R}^6$ and gives the nearly optimal lower bound $\Theta(\cC_{2,2})\geq \lambda(\mathbb{S}^2)$. Indeed, this estimate is sharp for any method that compares to the entropy of spheres because $\Theta(\cC_{2,2})<\lambda(\mathbb{S}^1)$ -- see Appendix \ref{Explicit}.    However,  it does not give any information about the density of $\cC_{1,3}\subset \Real^6$ even though direct computations give $\Theta(\cC_{1,3})>\lambda(\mathbb{S}^1)>\Theta(\cC_{2,2})$.  Moreover, even if one were able to strengthen the topological conclusions so that it gave information about $\cC_{1,3}$ the result would still not recover the bound $\lambda(\mathbb{S}^1)$ and so would not be sharp in that regard.  This is in contrast to Proposition \ref{HomolBallProp} which applies to all generalized Simons' cones.
\end{rem}

\begin{proof}[Proof of Proposition \ref{WeakProp}]
     Arguing as in the proof of Proposition \ref{HomolBallProp} gives that $\Omega_\mp(\Sigma_\pm)$ has trivial reduced homology groups and when $2\leq k\leq n-3$, $H_k(\Omega_\pm(\Sigma_\pm))$ is trivial. By Lemma \ref{StarShapeExpanderLem}, if $\omega_\pm=\mathrm{int}(\Omega_\pm(\mathcal{C}))\cap\mathbb{S}^n$, then $\Omega_\pm(\Sigma_\pm)$ and $\omega_\pm$ are homotopy equivalent and so have the same homology groups. Thus, combined with the Alexander duality (see \cite[Theorem 3.44]{HatAlgebraTopol}), it follows that, for $k\geq n-1$, all $H_k(\Omega_\pm(\Sigma_\pm))$ are trivial. Likewise, by the universal coefficient theorem -- e.g,  \cite[Theorem 3.2]{HatAlgebraTopol} -- both $H_1(\Omega_\pm(\Sigma_\pm))$ and $H_{n-2}(\Omega_\pm(\Sigma_\pm))$ are torsion-free. Hence, the claim on the homology of $\Omega_\pm(\Sigma_\pm)$ holds. To complete the proof, observe that by Lemma \ref{StarShapeExpanderLem} both $\Sigma_\pm$ and $\omega_\pm$ are homotopy equivalent to $\Omega_\pm(\Sigma_\pm)$ and so have the same homology groups as $\Omega_\pm(\Sigma_\pm)$. The description of the homology of $\mathcal{L}(\mathcal{C})$ follows immediately from this and the Mayer--Vietoris sequence (see \cite[p. 149]{HatAlgebraTopol}).
\end{proof}

\appendix

\section{Explicit Densities and Entropies} \label{Explicit}
We give the explicit value of the densities of the generalized Simons' cones as well as the entropies of spheres. We first remark that it follows from a computation of Stone \cite[Appnedix A]{StoEntSphere} that
$$
\lambda(\mathbb{S}^n)=2 \sqrt{\pi}\left(\frac{n}{2e}\right)^{\frac{n}{2}}\frac{1}{\Gamma\left( \frac{n+1}{2}\right)}=\left(\frac{n}{2\pi e} \right)^{\frac{n}{2}} \sigma_n
$$
where $\sigma_n$ is the area of $\mathbb{S}^n$ and is given by
$$
\sigma_n=(n+1) \omega_{n+1}=(n+1)\frac{\pi^{\frac{n+1}{2}}}{\Gamma(\frac{n+3}{2})}.
$$
Likewise, Ilmanen--White \cite{IWDenMinCone} compute that for $\cC_{m,l}\subset \Real^{m+l+2}=\Real^{n+1}$, 
$$
\Theta(\cC_{m,l})=\lambda(\cC_{m,l})=\frac{\sigma_m \sigma_n}{\sigma_{m+l}}\left( \frac{m}{m+l}\right)^{\frac{m}{2}} \left( \frac{l}{m+l}\right)^{\frac{l}{2}}.
$$
As observed in \cite{IWDenMinCone},
$$
\lim_{l\to \infty} \Theta(\cC_{m,l})=\lambda(\mathbb{S}^m),  \lim_{m\to \infty} =\Theta(\cC_{m,l})=\lambda(\mathbb{S}^l)
$$
and
$$
\lim_{m,l\to \infty} \Theta(\cC_{m,l})=\lim_{k\to \infty} \lambda(\mathbb{S}^k)=\sqrt{2}.
$$
The densities of generalized Simons' cones in low dimensions is recorded in Table \ref{TableDensities}.
\begin{table}[t]
	\caption{Densities of generalized Simons' cones}
		\label{TableDensities}
	\begin{tabular}{|c|c c c c c c|}
		\hline	$\Theta(\cC_{m,l})$ & $m=1$ & $m=2$ & $m=3$ & $m=4$ & $\cdots$ & $\begin{array}{c}\lim\limits_{m\to \infty} \Theta(\cC_{m,l})\\ =\lambda(\mathbb{S}^l)\end{array}$	\\
		\hline	$l=1$ &	$\approx 1.57$& $\approx 1.54$  &
		$ \approx 1.530$  & $\approx 1.526$ &$\cdots$ &  $\approx 1.52$   \\
		$l=2$ &		$\approx 1.54$	& $1.5$ &  $\approx 1.487$  & $\approx 1.481 $&  $\cdots$ & $\approx 1.47$ \\
		$l=3$ &		$\approx 1.530$	& $\approx 1.487$ & $\approx 1.473 $ &  $\approx 1.466$  & $\cdots$ & $\approx 1.453$ \\
		$l=4$ &		$\approx 1.526$	& $\approx 1.481 $ & $\approx 1.466$ &  $\frac{35}{24}\approx 1.46$  & $\cdots$  &$\approx 1.444$ \\
		$\vdots$ &  $\vdots$ &  $\vdots$   &  $\vdots$ &  $\vdots$ & $\ddots$  & $\vdots$ \\
		$\begin{array}{c}\lim\limits_{l\to \infty} \Theta(\cC_{m,l})\\ =\lambda(\mathbb{S}^m)\end{array}$	 &$\approx 1.52$ & $\approx 1.47$ & $\approx 1.453$ & $\approx 1.444$ & $\cdots$  & $\sqrt{2}\approx 1.414$\\
		\hline
	\end{tabular}

\end{table}

It is also convenient to set
$$
Z(n)=\frac{\lambda(\mathbb{S}^{n-2})}{\lambda(\mathbb{S}^{n-1})}
$$  
which is the density lower bounds proved by Zhu in \cite{ZhuThesis}. We compare the densities of generalized Simons' cones, $\cC_{k, n-k-1}\subset \Real^{n+1}$ for $3\leq n \leq 7$ and $1\leq k \leq \lfloor\frac{n}{2}\rfloor$ with the bounds provided by Propositions \ref{ContractProp}, \ref{HomolBallProp} and \ref{WeakProp} and $Z(n)$ in Table \ref{BoundTable}.

\begin{table}[ht]
	\caption{ Densities of generalized Simons' cones and theoretical bounds}
	\label{BoundTable}
\begin{tabular}{|c|c c c c c c|}\hline
$\Theta(\cC_{k,n-k-1})$ & $n=3$ & $n=4$ & $n=5$ & $n=6$  & $\cdots$ & $n\to \infty$ \\
\hline 
$k=1$ &	$\approx 1.57$& $\approx 1.54$   & $ \approx 1.530$  & $\approx 1.526$ &  $\cdots$ & $\approx 1.52$ \\
 $k=2$ & * & * & $ 1.5$ &  $\approx 1.487$  &$\cdots$ & $\approx 1.47$\\
 $\vdots$ & *  & * & * & * &  $\ddots$ & $\vdots $\\
 $k\to \infty$ & * & * & * & * & $\cdots$ & $\sqrt{2}\approx 1.414$ \\ \hline
 $\lambda(\mathbb{S}^{n-3}) $& * & * & $\approx 1.47$ & $\approx 1.453$  &$\cdots$ & $\sqrt{2}\approx 1.414$ \\

$\lambda(\mathbb{S}^{n-2})$  &$\approx 1.52$ & $\approx 1.47$ & $\approx 1.453$ & $\approx 1.444$ &$\cdots$ & $\sqrt{2}\approx 1.414$ \\

$\lambda(\mathbb{S}^{n-1})$ & $\approx 1.47$ & $\approx 1.453$ & $\approx 1.444$ & $\approx 1.438$ &$\cdots$ & $\sqrt{2}\approx 1.414$ \\
 
$Z(n)$ & $\approx 1.03$& $\approx 1.01$ & $\approx 1.007$ & $\approx 1.004$  &$\cdots$ & $1$ \\
\hline
\end{tabular}
\end{table}

Finally, we compare the bounds using entropy of spheres with the universal bounds provided by \cite{ChengLiYau}.  Let
$$
C_n'=\frac{n^{\frac{n}{2}} e \Gamma\left(\frac{n}{2},1\right)}{2}
$$
where $\Gamma(s,x)$ is the incomplete Gamma function and consider the explicit constants
$$
\epsilon_{CLY}(n)=\frac{1}{2n+3+2 \exp(2n C_n')}>0.
$$
Cheng--Li--Yau \cite{ChengLiYau} show that if $\cC\subset \Real^{n+1}$ is a non-flat regular minimal cone, then
$$
\Theta(\cC)\geq 1+\epsilon_{CLY}(n).
$$
These numbers are very small.  Indeed,  for $n\geq 3$, 
$$
 2.3\times 10^{-10}\approx \epsilon_{CLY}(3)> \epsilon_{CLY}(n). 
$$

\end{document}